\documentclass[12pt]{amsart}
\usepackage{amsmath, amscd, amssymb, graphicx}

\numberwithin{equation}{section}
\newtheorem{theorem}{Theorem}[section]

\newtheorem{conjecture}[theorem]{Conjecture}
\theoremstyle{remark}

\theoremstyle{definition}

\begin{document}
\title[New structure for orthogonal quantum group invariants]{New structure for orthogonal quantum group invariants}

\author{Qingtao Chen}
\address{Q. Chen, Mathematics Section, International Center for Theoretical Physics, Strada Costiera, 11, Trieste, I-34151, Italy (qchen1@ictp.it)}

\author{Kefeng Liu}
\address{K. Liu, Center of Mathematical Sciences, Zhejiang University, Box 310027, Hangzhou, P.R.China. Department of mathematics, University of California at Los Angeles, Box 951555\newline
Los Angeles, CA 90095-1555, USA (liu@math.ucla.edu)}

\date{}

\begin{abstract}
Based on the orthogonal Labastida-Mari{\~n}o-Ooguri-Vafa conjecture made by L. Chen \& Q. Chen \cite{CC}%
, we derive an infinite product formula for Chern-Simons partition functions, which generalizes the Liu-Peng's \cite{LP3} recent results to the orthogonal case. Symmetry property of this new infinite product structure is also discussed.

\end{abstract}

\maketitle

\section{Introduction}

\label{sec1}In 1984 Jones \cite{J1, J2} discovered a polynomial invariant of
oriented knots. Later, along this clue, the two-variable HOMFLY-PT \cite%
{HOMFLY, PT} and one-variable Kauffman polynomials \cite{Kau1} of oriented
and unoriented links have also been discovered in 1985 and 1987
respectively. HOMFLY-PT polynomial generalizes both Alexander and Jones
polynomial. In 1990, two-variable Kauffman polynomial \cite{Kau2} was also
introduced, which only generalizes the Jones polynomial.

By using the Chern-Simons path integral method Witten \cite{W} gave a
quantum field theory interpretation of the Jones polynomial. Witten \cite{W}
also predicted the existence of 3-manifold quantum invariants. Following
this track, Reshetikhin and Turaev \cite{RT1, RT2} gave a construction of
3-manifold invariants by using quantum universal enveloping algebra (quantum
group) $U_{q}(sl_{2})$ at roots of unity, which leads to the colored version
of classical HOMFLY-PT and Kauffman polynomial invariants. These work
actually give a unified understanding of the quantum group invariants of
links. Here color means the representation of quantum groups. It turns out
that colored HOMFLY-PT polynomial is a special linear quantum group
invariants, i.e. the quantum group of $A_{n}$ type; while colored Kauffman
is a quantum group invariants of $B_{n}$, $C_{n}$ and $D_{n}$ type.

In a series of papers, Labastida, Mari\~{n}o, Ooguri and Vafa \cite{LaM,
LMV, OV} proposed a conjectural description of deep relationship between
reformulated invariants of colored HOMFLY-PT link in late 1990s and early
2000s. This conjecture was proved by K. Liu and P. Peng in \cite{LP1, LP2}.

In some sense the Labastida-Mari\~{n}o-Ooguri-Vafa (LMOV) conjecture can be
expressed purely by using mathematical language, i.e. irreducible
representation of quantum groups. The physics background of this conjecture
can be dated back to 't Hooft's seminal work on large $N$ expansion of $U(N)$
gauge field theories in 1974. Gopakumar and Vafa \cite{GV} described the
exact theory that closed topological string theory on the resolved conifold
is dual to the $U(N)$ Chern-Simons theory on $S^{3}$. The Gromov-Witten
theory of the resolved conifold actually corresponds to the Chern-Simons
theory of an unknot. The LMOV conjecture considers the more general case
when the link or knot is nontrivial and the corresponding Wilson loop
expectation values, i.e. colored HOMFLY-PT polynomial of the link. So the
LMOV conjecture could be viewed as a counterpart of Gopakumar-Vafa
conjecture.

Previously people thought that only HOMFLY-PT can be expressed as a series
in $q-q^{-1}$ and $t$ with integer coefficients; while the colored HOMFLY-PT
polynomial is just a polynomial of $q^{\pm 1}$ and $t^{\pm 1}$ with rational
coefficients. The LMOV conjecture predicts an intrinsic symmetry of $%
q-q^{-1} $ about reformulated invariants of the colored HOMFLY-PT polynomial
and hidden integrality encoded in the colored HOMFLY-PT polynomial.

In the late 2010s, the orthogonal LMOV conjecture was also formulated by Lin
Chen \& Qingtao Chen \cite{Che1, CC, Che2}. and Marcos Mari\~{n}o \cite{M}.
Chen-Chen's formulation put an emphasis on the colored Kauffman solely;
while Mari\~{n}o's put an emphasis on the relation between the colored
HOMFLY-PT polynomials and the colored Kauffman polynomials.

More recently, Kefeng Liu and Pan Peng \cite{LP3} obtained a new structure
of the colored HOMFLY-PT polynomial that the Chern-Simons partition function
appeared in the LMOV conjecture can be expressed as an infinite product,
which indicates some potential modularity of the Chern-Simons partition
function \cite{LP3}.

In this paper, an infinite product expression for the orthogonal
Chern-Simons partition function appeared in the orthogonal LMOV type
conjecture is established and the case of unknot is presented in explicit
formula.

This paper is organized as follows. In Section 2, we introduce the basic
setups for the quantum group invariants of links. In Section 3, we describe
the original and the orthogonal LMOV conjecture and notations. In Section 4,
we derive the orthogonal Chern-Simons partition function as an infinite
product and illustrate an example of the unknot. In Section 5, we discuss
the symmetric properties associated to this infinite product structure.

\section{Quantum Invariants of Links}

\label{sec2}Let $\mathfrak{g}$ be a finite dimensional complex semi-simple
Lie algebra of rank $N$ with Cartan matrix $(C_{ij})$. Let $U_{q}(\mathfrak{g%
})$ be the quantum enveloping algebra of $\mathfrak{g}$. Let $V$ be a vector
space over a field $k$. A linear automorphism $c$ of $V\otimes V$ is said to
be an $R$-matrix if it is a solution of the following Yang-Baxter equation

\begin{equation}
(c\otimes id_{V})(id_{V}\otimes c)(c\otimes id_{V})=(id_{V}\otimes
c)(c\otimes id_{V})(id_{V}\otimes c)
\end{equation}%
that holds in the automorphism group of $V\otimes V\otimes V$.

It is well known that the solution of the Yang-Baxter equation provides a
representation of the braid group. The solution we used is the following
so-called universal $R$-matrix.

\begin{equation}
R=q^{\underset{i,j}{\sum }C_{ij}^{-1}H_{i}\otimes H_{j}}\underset{\beta }{%
\prod }exp_{q}[(1-q^{-2})X_{\beta }^{+}\otimes X_{\beta }^{-}],
\end{equation}%
where $\beta $ runs over positive roots of $sl(N,\mathbf{C})$, $(C_{ij})$ is
the Cartan matrix, and $q-$exponential is given by

\begin{equation}
exp_{q}[x]=\underset{k=0}{\overset{\infty }{\sum }}q^{\frac{1}{2}k(k+1)}%
\frac{x^{k}}{[k]_{q}!},
\end{equation}%
where $[k]_{q}!=[k]_{q}\cdot \lbrack k-1]_{q}\cdot \cdot \cdot \lbrack
1]_{q} $, $[k]_{q}=\dfrac{[k]}{[1]}$ and $[n]=q^{n}-q^{-n}$.

Given a link $\mathcal{L}$ with $L$ components. It is well-known that $%
\mathcal{L}$ can be represented by an element in some braid group $B_{m}$
with $m$ strands. For each component, we associate to it with an irreducible
representation $A^{\alpha }$ of quantized universal enveloping algebra $%
U_{q}(sl(N,\mathbb{C}))$. $A^{\alpha }$ is labeled by the highest weight $%
\Lambda _{\alpha }$. As usual, we associate to them with the Young diagrams.
Without loss of generality, one can assume the first $m_{1}$ strands
correspond to the first component, the second $m_{2}$ strands correspond to
the second component, and so on. Let

\begin{equation}
\widehat{V}=\overset{L}{\underset{\alpha =1}{\bigotimes }}V_{\Lambda
_{\alpha }}^{m_{\alpha }}.
\end{equation}
And write the braiding $\check{R}=P_{12}R:\, V\otimes W\rightarrow W\otimes
V $, where $P_{12}(s\otimes t)=t\otimes s$.

For a generator of the braid group $B_{m}$, $\sigma _{i}$, define

\begin{equation}
\pi (\sigma _{i}^{\pm 1})=Id_{V_{1}}\otimes Id_{V_{2}}\otimes \cdot \cdot
\cdot \otimes Id_{V_{i-1}}\otimes \check{R}^{\pm 1}\otimes \cdot \cdot \cdot
\otimes Id_{V_{m}}.
\end{equation}%
The quantum group invariants of the link $\mathcal{L}$ is defined as follows

\begin{equation}
W_{A^{1},\cdot \cdot \cdot ,A^{L}}^{\mathfrak{g}}(\mathcal{L})=q^{d(\mathcal{%
L})}Tr_{\widehat{V}}(\mu ^{m}\cdot \pi (\mathcal{L})),
\end{equation}%
where $\mu =q^{\rho ^{\ast }}$, $\rho ^{\ast }$ is the element in $\mathfrak{%
h\subset }U_{q}(\mathfrak{h})$ corresponding to Weyl vector (i.e. the sum of
fundamental weights) under the natural isomorphism $\mathfrak{h\simeq h}%
^{\ast }$ and $d(\mathcal{L})$ is given by the following formula

\begin{equation}
d(\mathcal{L})=-\underset{\alpha =1}{\overset{L}{\sum }}\omega (\mathcal{K}%
_{\alpha })(\Lambda _{\alpha },\Lambda _{\alpha }+2\rho )+\frac{2}{N}%
\underset{\alpha <\beta }{\overset{L}{\sum }}lk(\mathcal{K}_{\alpha },%
\mathcal{K}_{\beta })l_{\alpha }l_{\beta }.
\end{equation}

Special case 1: For the unknot $\bigcirc $, $W_{A}(\bigcirc )$ is the
quantum dimension $\dim _{q}(V_{A})$ of the corresponding representation
space $V_{A}$.

Special case\ 2: If $\mathfrak{g}=sl_{N}$ and $A^{1}=A^{2}=\cdot \cdot \cdot
=A^{L}=(1)$, the quantum group invariant of links equal to the HOMFLY
polynomial at $t=q^{N}$ up to a universal factor $\dfrac{t-t^{-1}}{q-q^{-1}}$%
.

Special case\ 3: If $\mathfrak{g}=so_{2N+1}$ and $A^{1}=A^{2}=\cdot \cdot
\cdot =A^{L}=(1)$, quantum group invariant of links equal to Kauffman
polynomial at $t=q^{2N}$ up to a universal factor $1+\dfrac{t-t^{-1}}{%
q-q^{-1}}$ and some $t$ power of the linking numbers.

Thus the quantum group invariant associated to $\mathfrak{g}=sl_{N}$ and $%
\mathfrak{g}=so_{2N+1}$ are called the colored HOMFLY and the colored
Kauffman polynomials respectively.

Actually the irreducible representation of the quantum groups of special
linear and orthogonal cases can be labeled by the Young Tableau. Now we
would like to introduce some basic notation of the partition and the
corresponding Young Tableau.

A partition of $n$ is a tuple of positive integers $\mu =(\mu _{1},\mu
_{2},\ldots ,\mu _{k})$ such that $|\mu |\triangleq \overset{k}{\underset{i=1%
}{\sum }}\mu _{i}=n$ and $\mu _{1}\geq \mu _{2}\geq \cdots \geq \mu _{k}>0$,
where $|\mu |$ is called the degree of $\mu $ and $k$ is called the length
of $\mu $, denoted by $\ell (\mu )$. A partition can be represented by a
Young diagram, for example, partition $(5,4,2,1)$ can be identified as the
following Young diagram.

Denote by $\mathcal{P}$ the set of all Young diagrams. Let $\chi _{A}$ be
the character of irreducible representation of symmetric group, labeled by
partition $A$ . Given a partition $\mu $, define $m_{j}=\#(\mu _{k}=j;k\geq
1)$. The order of the conjugate class of type $\mu $ is given by:

\begin{equation}
\mathfrak{z}_{\mu }=\prod_{j\geq 1}j^{m_{j}}m_{j}!.
\end{equation}
The theory of symmetric functions has a close relationship with the
representations of symmetric group. The symmetric power functions of a given
set of variables $x=\{x_{j}\}_{j\geq 1}$ are defined as the direct limit of
the Newton polynomials:

\begin{equation}
p_{n}(x)=\overset{\infty }{\sum_{j=1}}x_{j}^{n},\qquad p_{\mu }(x)=\overset{%
\ell (\mu )}{\prod_{i=1}}p_{\mu _{i}}(x).
\end{equation}

We will consistently denote by $\mathcal{L}$ a link and by $L$ the number of
components in $\mathcal{L}$.

The irreducible $U_{q}(\mathfrak{g})$ modules associated to $\mathcal{L}$
will be labeled by their highest weights, thus by Young diagrams. We usually
denote it by a vector form $\vec{A}=(A^{1},\ldots ,A^{L})$.

Let $\vec{x}=(x^{1},\ldots ,x^{L})$ is $L$ sets of variables, each of which
is associated to a component of $\mathcal{L}$ and $\vec{\mu}=(\mu
^{1},\ldots ,\mu ^{L})\in \mathcal{P}^{L}$ be a tuple of $L$ partitions.
Define:
\begin{align*}
& [\mu ]=\prod_{i=1}^{\ell (\mu )}[\mu _{i}], & & [\vec{\mu}]=\prod_{\alpha
=1}^{L}[\mu ^{\alpha }], & & \mathfrak{z}_{\vec{\mu}}=\prod_{\alpha =1}^{L}%
\mathfrak{z}_{\mu ^{\alpha }}, \\
& \chi _{\vec{A}}(C_{\vec{\mu}})=\prod_{\alpha =1}^{L}\chi _{A^{\alpha
}}(C_{\mu ^{\alpha }}), & & s_{\vec{A}}(\vec{x})=\prod_{\alpha
=1}^{L}s_{A^{\alpha }}(x^{\alpha }), & & p_{\vec{\mu}}(\vec{x}%
)=\prod_{\alpha =1}^{L}p_{\mu ^{\alpha }}(x^{\alpha }).
\end{align*}

When we consider the orthogonal quantum group invariants, we need to study
the Brauer algebra $Br_{n}$ which contains the group algebra $\mathbb{C}%
[S_{n}]$ as a direct summand. Thus all the irreducible representations of $%
S_{n}$ are also irreducible representations of $Br_{n}$, labeled by
partitions of the integer $n$. Indeed, the set of irreducible
representations of $Br_{n}$ are bijective to the set of partitions of the
integers $n-2k$, where $k=0,1,\cdots ,[\frac{n}{2}]$ \cite{Ram2, Wen}. Thus
the semi-simple algebra $Br_{n}$ can be decomposed into the direct sum of
simple algebras

\begin{equation}
Br_{n}\cong \overset{\lbrack \frac{n}{2}]}{\bigoplus_{k=0}}%
\bigoplus_{\lambda \vdash n-2k}M_{d_{\lambda }\times d_{\lambda }}(\mathbb{C}%
).
\end{equation}

The work of Beliakova and Blanchet \cite{BB} constructed an explicit basis
of the above decomposition. An up and down tableau $\Lambda =(\lambda
_{1},\lambda _{2},\cdots ,\lambda _{n})$ is a tube of $n$ Young diagrams
such that $\lambda _{1}=(1)$ and each $\lambda _{i}$ is obtained by adding
or removing one box from $\lambda _{i-1}$. Let $\lambda $ be a partition of $%
n-2k$. Denote by $|\Lambda |=\lambda $ if $\lambda _{n}=\lambda $, and we
say an up and down tableau $\Lambda $ is of shape $\lambda $. There is a
minimal path idempotent $p_{\Lambda }\in Br_{n}$ associated to each $\Lambda
$. Then the minimal central idempotent $\pi _{\lambda }$ of $Br_{n}$
correspond to the irreducible representation labeled by $\lambda $ is given
by

\begin{equation}
\pi _{\lambda }=\sum_{|\Lambda |=\lambda }p_{\Lambda }.
\end{equation}
In particular, the dimension of the irreducible representations $d_{\lambda
} $ is the number of up and down tableau of shape $\lambda $. More detail
can be found in \cite{BB, Wen}.

The table of the characters and orthogonal relations can be found in \cite%
{Ram1, Ram2, Ram3}. The values of a character of $Br_{n}$ is completely
determined by its values on the set of elements $e^{k}\otimes \gamma
_{\lambda }$, where $e$ is the conjugacy class of $e_{1},\cdots ,e_{n-1}$
and $\gamma _{\lambda }$ is the conjugacy class in $S_{n-2k}$ labeled by the
partition $\lambda $ of $n-2k$. The notion $e^{k}\otimes \gamma _{\lambda }$
stands for the tangle in the following diagram.
\begin{align*}
& \ e_{0}\ \ \ \ \ e_{2}\ \ \ \ \cdots \ \ \ \ \ e_{2k}\ \ \ \ \ \ \gamma
_{\lambda }\ \ \  \\
& {%
\includegraphics[height=1.0153in, width=2.1318in]
{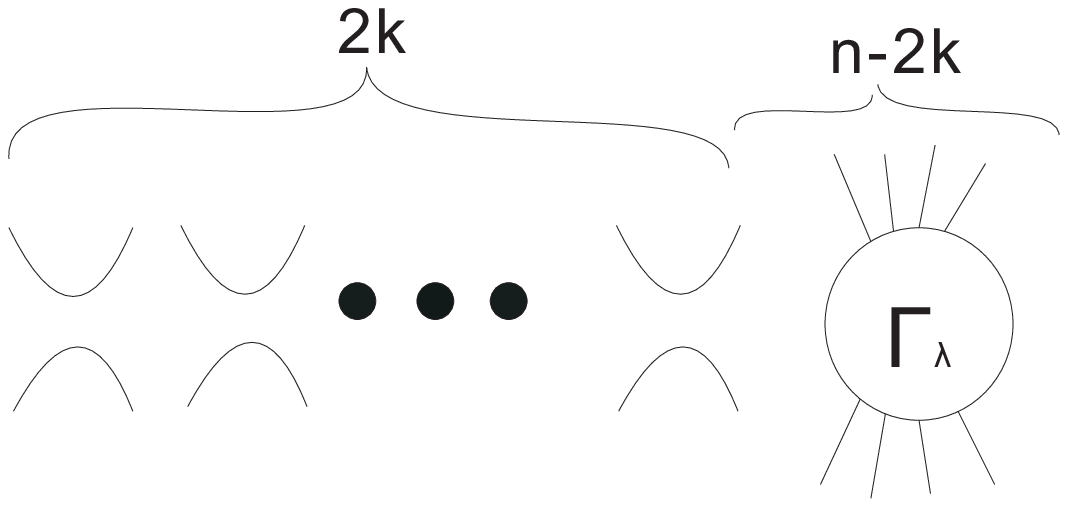}}
\end{align*}%
where $\Gamma _{\lambda }$ is a diagram in the conjugacy class of $S_{n-2k}$
labeled by a partition $\lambda $ of $n-2k$.

Denote by $\chi _{A}$ the character of the irreducible representation of $%
Br_{n} $ labeled by a partition $A\vdash n-2k$ for some $k$, and denote by $%
\chi _{B}^{S_{n}}$ the character of the irreducible representation of $S_{n}$
labeled by a partition $B\vdash n$. It is known that when $A$ is a partition
of $n$, then $\chi _{A}(e^{m}\otimes \gamma _{\lambda })=0$ for all $m>0$
and partition $\lambda \vdash n-2m$, and $\chi _{A}(\gamma _{\mu })=\chi
_{A}^{S_{n}}(\gamma _{\mu })$ for partition $\mu \vdash n$ coincide with the
characters of the permutation group $S_{n}$ \cite{Ram2}.

\section{Labastida-Mari\~{n}o-Ooguri-Vafa type Conjecture for colored
Kauffman polynomial}

\label{sec3}Let's quickly review the original LMOV conjecture first.

For each link $\mathcal{L}$, the type$-A$ Chern-Simons partition function of
$\mathcal{L}$ is defined by

\begin{equation}
Z_{CS}^{SL}(\mathcal{L};q,t;\overrightarrow{x})=\sum_{\overrightarrow{A}\in
\mathcal{P}^{L}}W_{\overrightarrow{A}}^{SL}(\mathcal{L};q,t)s_{%
\overrightarrow{A}}(\overrightarrow{x})=1+\sum_{\overrightarrow{\mu }\neq
\overrightarrow{0}}Z_{\overrightarrow{\mu }}^{SL}p_{\overrightarrow{\mu }}(%
\overrightarrow{x}),
\end{equation}%
where $s_{\overrightarrow{A}}(\overrightarrow{x})$ are the Schur polynomials.

The original LMOV conjecture describes a very subtle structure of $%
Z_{CS}^{SL}(\mathcal{L};q,t;\overrightarrow{x})$, which was proved by Kefeng
Liu and Pan Peng \cite{LP1, LP2}, based on the cabling technique and a
careful degree analysis of the cut-join equation. As an application, the
LMOV conjecture gives highly non-trivial relations between colored HOMFLY
polynomials. The first such relation is the classical Lichorish-Millett
theorem \cite{LiM}.

The study of the colored Kauffman polynomials are more difficult. For
instance, the definition of the Chern-Simons partition function for the
orthogonal quantum groups involves the representations of Brauer centralizer
algebras, which admit a more complicated orthogonal relations \cite{Ram1,
Ram2, Ram3}. In a joint work with Lin Chen, we \cite{CC} rigorously
formulate the orthogonal quantum group version of LMOV conjecture by using
the representation of the Brauer centralizer algebra.

Now we set $\mathfrak{g=}so_{2N+1}$.

Let $\mathcal{L}$ be a link with $L$ components and%
\begin{equation}
pb_{n}(z)=\overset{\infty }{\underset{j=-\infty }{\sum }}z_{j}^{n},\qquad
pb_{\mu }(z)=\overset{\ell (\mu )}{\prod_{i=1}}pb_{\mu _{i}}(z),\qquad pb_{%
\overrightarrow{\mu }}(\overrightarrow{z})=\overset{L}{\prod_{\alpha =1}}%
pb_{\mu ^{\alpha }}(z^{\alpha }).
\end{equation}

Let $Z_{CS}^{SO}(\mathcal{L},q,t)$ be the orthogonal Chern-Simons partition
function defined by

\begin{equation}
Z_{CS}^{SO}(\mathcal{L};q,t;\overrightarrow{z})=\sum_{\overrightarrow{\mu }%
\in \mathcal{P}^{L}}\frac{pb_{\overrightarrow{\mu }}(\overrightarrow{z})}{%
\mathfrak{z}_{\vec{\mu}}}\underset{\overrightarrow{A}\in \widehat{Br}_{|%
\overrightarrow{\mu }|}}{\sum }\chi _{\overrightarrow{A}}(\gamma _{%
\overrightarrow{\mu }})W_{\overrightarrow{A}}^{SO}(\mathcal{L};q,t),
\end{equation}%
where $\mathfrak{z}_{\vec{\mu}}=\frac{||\overrightarrow{\mu }||!}{|C_{%
\overrightarrow{\mu }}|}$, $|\overrightarrow{\mu }|=(d^{1},...,d^{L})$, $%
\widehat{Br}_{|\overrightarrow{\mu }|}$ denotes the set $\widehat{Br}%
_{d^{1}}\times \cdot \cdot \cdot \times \widehat{Br}_{d^{L}}$ (every element
is a representation of the Brauer algebra), $\overrightarrow{\mu }=(\mu
^{1},...,\mu ^{L})$ for partitions $\mu ^{i}$ of $d^{i}\in \mathbb{Z}$ and $%
\chi _{\overrightarrow{A}}(\gamma _{\overrightarrow{\mu }})=\overset{L}{%
\underset{i=1}{\prod }}\chi _{A^{i}}(\gamma _{\mu ^{i}})$ for the character $%
\chi _{A^{i}}$ of $Br_{d^{i}}$ labeled by $A^{i}$.

Expend the free energy

\begin{equation}
F_{CS}^{SO}(\mathcal{L};q,t;\overrightarrow{z})=\log Z_{CS}^{SO}(\mathcal{L}%
;q,t;\overrightarrow{z})=\sum_{\overrightarrow{\mu }\neq \overrightarrow{0}%
}F_{\overrightarrow{\mu }}^{SO}(\mathcal{L};q,t)pb_{\overrightarrow{\mu }}(%
\overrightarrow{z}),
\end{equation}

Then the reformulated invariants are defined by

\begin{equation}
g_{\overrightarrow{\mu }}(\mathcal{L};q,t)=\sum_{k|\overrightarrow{\mu }}%
\frac{\mu (k)}{k}F_{\overrightarrow{\mu }/k}(\mathcal{L};q^{k},t^{k}).
\end{equation}

The orthogonal LMOV conjecture was formulated by L. Chen and Q. Chen \cite%
{CC} as follows

\begin{conjecture}[Orthorgonal LMOV, Chen-Chen \protect\cite{CC}]
\label{Main}%
\begin{equation*}
\frac{\mathfrak{z}_{\vec{\mu}}[1]^{2}\cdot \lbrack g_{\overrightarrow{\mu }}(%
\mathcal{L};q,t)-g_{\overrightarrow{\mu }}(\mathcal{L};q,-t)]}{2[%
\overrightarrow{\mu }]}=\overset{\infty }{\sum_{g=0}}\sum_{\beta \in \mathbb{%
Z}}N_{\overrightarrow{\mu },g,\beta }(q-q^{-1})^{g}t^{\beta },
\end{equation*}%
where $N_{\overrightarrow{\mu },g,\beta }$ are the integer coefficients and
vanish for sufficiently large $g$ and $|\beta |$,
\end{conjecture}

This conjecture is a rigorous mathematical formulation of the LMOV type
conjecture about the colored Kauffman polynomial; while in \cite{BFM, M},
their conjecture emphasizes on the relationship between colored HOMFLY-PT
and colored Kauffman. The integer coefficients $N_{\overrightarrow{\mu }%
,g,\beta }$ are closely related to the BPS numbers.

\section{Infinite product formula for orthogonal Chern-Simons partition
functions}

\label{sec4}To derive an infinite product formula, we will state the result
for a knot first, since the notations in the computation for a knot are
relatively simpler.

\subsection{The case of a knot}

Given $z=\{z_{i}\}_{-\infty <i<\infty }$, $x=\{x_{j}\}_{j\geq 1}$, define
\begin{equation*}
x\ast y=\{x_{i}\cdot y_{j}\}_{-\infty <i<\infty ,j\geq 1}.
\end{equation*}%
We also define $z^{d}=\{z_{i}^{d}\}_{-\infty <i<\infty }$. The $d$-th Adam
operation of a type-B Schur function is given by $sb_{A}(z^{d})$.

Introduce variables $q^{\rho }=\{-q^{2j-1}\}_{j\geq 1}$. we have

\begin{equation}
p_{n}(q^{\rho })=\frac{1}{[n]},
\end{equation}%
where we assume $|q|<1$.

Set $w=z\ast q^{\rho }$, then we have

\begin{equation}
pb_{n}(w)=\frac{pb_{n}(z)}{[n]}
\end{equation}%
and

\begin{equation}
pb_{\mu }(w)=\frac{pb_{\mu }(z)}{[\mu ]}.
\end{equation}

Consider the free energy weighted by the $pb_{\mu }(w)$, orthogonal LMOV
conjecture implies the following reformulation of the free energy:

\begin{equation}
F^{SO}(\mathcal{K};q,t;w)=\overset{\infty }{\underset{d=1}{\sum }}\underset{%
\mu \neq 0}{\sum }\frac{1}{d}g_{\mu }(\mathcal{K};q^{d},t^{d})pb_{\mu
}(w^{d})
\end{equation}

and

\begin{equation}
g_{\mu }(\mathcal{K};q,t)-g_{\mu }(\mathcal{K};q,-t)=\frac{2[\mu ]}{%
\mathfrak{z}_{\mu }[1]^{2}}\overset{\infty }{\sum_{g=0}}\sum_{\beta \in
\mathbb{Z}}N_{\overrightarrow{\mu },g,\beta }(q-q^{-1})^{g}t^{\beta }.
\end{equation}

There exists integers $n_{B;\,g,\beta }$ such that

\begin{equation}
\overset{\infty }{\sum_{g=0}}N_{B,g,\beta }(q-q^{-1})^{g}=\overset{\infty }{%
\sum_{g=0}}n_{B,g,\beta }\underset{k=0}{\overset{g}{\sum }}q^{g-2k}.
\end{equation}

By the orthogonal LMOV conjecture, $N_{B;g,\beta }$ vanish for sufficiently
large $g$ and $|\beta |$, thus $n_{B;\,g,\beta }$ vanish for sufficiently
large $g$ and $|\beta |$.

We have

\begin{eqnarray}
&&F^{SO}(\mathcal{K};q,t;w)-F^{SO}(\mathcal{K};q,-t;w)  \notag \\
&=&\overset{\infty }{\underset{d=1}{\sum }}\underset{\mu \neq 0}{\sum }\frac{%
1}{d}g_{\mu }(\mathcal{K};q^{d},t^{d})pb_{\mu }(w^{d})-\overset{\infty }{%
\underset{d=1}{\sum }}\underset{\mu \neq 0}{\sum }\frac{1}{d}g_{\mu }(%
\mathcal{K};q^{d},(-1)^{d}t^{d})pb_{\mu }(w^{d})  \notag \\
&=&\underset{d\in O\mathbb{Z}_{+}}{\sum }\underset{\mu \neq 0}{\sum }\frac{1%
}{d}(g_{\mu }(\mathcal{K};q^{d},t^{d})-g_{\mu }(\mathcal{K}%
;q^{d},-t^{d}))pb_{d\mu }(w)  \notag \\
&=&\underset{d\in O\mathbb{Z}_{+}}{\sum }\underset{\mu \neq 0}{\sum }\frac{1%
}{d}\frac{2[d\mu ]}{\mathfrak{z}_{\mu }[d]^{2}}\overset{\infty }{\sum_{g=0}}%
\underset{\beta \in \mathbb{Z}}{\sum }N_{\mu ,g,\beta
}(q^{d}-q^{-d})^{g}t^{d\beta }pb_{d\mu }(w)  \notag \\
&=&\underset{d\in O\mathbb{Z}_{+}}{\sum }\underset{\mu \neq 0}{\sum }\frac{1%
}{d}\frac{2}{\mathfrak{z}_{\mu }[d]^{2}}\overset{\infty }{\sum_{g=0}}%
\underset{\beta \in \mathbb{Z}}{\sum }N_{\mu ,g,\beta
}(q^{d}-q^{-d})^{g}t^{d\beta }pb_{d\mu }(z)  \label{integer N} \\
&=&\underset{d\in O\mathbb{Z}_{+}}{\sum }\underset{\mu \neq 0}{\sum }\frac{1%
}{d}\frac{2}{\mathfrak{z}_{\mu }[d]^{2}}\overset{\infty }{\sum_{g=0}}%
\underset{\beta \in \mathbb{Z}}{\sum }\underset{k=0}{\overset{g}{\sum }}%
n_{\mu ,g,\beta }q^{(g-2k)d}t^{d\beta }pb_{d\mu }(z)  \notag \\
&=&\underset{d\in O\mathbb{Z}_{+}}{\sum }\underset{\mu \neq 0}{\sum }\frac{1%
}{d}\frac{2}{\mathfrak{z}_{\mu }}\underset{m=1}{\overset{\infty }{\sum }}%
mq^{2md}\overset{\infty }{\sum_{g=0}}\underset{\beta \in \mathbb{Z}}{\sum }%
\underset{k=0}{\overset{g}{\sum }}n_{\mu ,g,\beta }q^{(g-2k)d}t^{d\beta
}pb_{d\mu }(z)  \notag \\
&=&\underset{d\in O\mathbb{Z}_{+}}{\sum }\underset{\mu \neq 0}{\sum }\frac{1%
}{d}\frac{2}{\mathfrak{z}_{\mu }}\overset{\infty }{\sum_{g=0}}\underset{%
\beta \in \mathbb{Z}}{\sum }\underset{m=1}{\overset{\infty }{\sum }}\underset%
{k=0}{\overset{g}{\sum }}mn_{\mu ,g,\beta }q^{(g-2k+2m)d}t^{d\beta }pb_{d\mu
}(z)  \notag \\
&=&\underset{\mu \neq 0}{\sum }\overset{\infty }{\sum_{g=0}}\underset{\beta
\in \mathbb{Z}}{\sum }\underset{m=1}{\overset{\infty }{\sum }}\underset{k=0}{%
\overset{g}{\sum }}\frac{mn_{\mu ,g,\beta }}{\mathfrak{z}_{\mu }}\underset{%
d\in O\mathbb{Z}_{+}}{\sum }\frac{2}{d}q^{(g-2k+2m)d}t^{d\beta }pb_{d\mu
}(z),
\end{eqnarray}

where $O\mathbb{Z}_{+}=\{1,3,5,...\}$ denotes the set of all positive odd
integers.

Now we analyze the following computation in detail,

\begin{align*}
& \underset{d\in O\mathbb{Z}_{+}}{\sum }\frac{2}{d}q^{(g-2k+2m)d}t^{d\beta
}pb_{d\mu }(z) \\
& =\underset{d\in O\mathbb{Z}_{+}}{\sum }\frac{2}{d}q^{(g-2k+2m)d}t^{d\beta }%
\underset{j=1}{\overset{\ell (\mu )}{\prod }}\left( \underset{i=-\infty }{%
\overset{\infty }{\sum }}(z_{i})^{d\mu _{j}}\right)  \\
& =\underset{d\in O\mathbb{Z}_{+}}{\sum }\frac{2}{d}q^{(g-2k+2m)d}t^{d\beta }%
\underset{i_{1},...,i_{\ell (\mu )}}{\sum }(z_{i_{1}}^{\mu
_{1}}z_{i_{2}}^{\mu _{2}}\cdot \cdot \cdot z_{i_{\ell (\mu )}}^{\mu _{\ell
(\mu )}})^{d} \\
& =\underset{i_{1},...,i_{\ell (\mu )}}{\sum }\underset{d\in O\mathbb{Z}_{+}}%
{\sum }\frac{2}{d}(q^{g-2k+2m}t^{\beta }z_{i_{1}}^{\mu _{1}}z_{i_{2}}^{\mu
_{2}}\cdot \cdot \cdot z_{i_{\ell (\mu )}}^{\mu _{\ell (\mu )}})^{d}
\end{align*}

Now we compute the series $\underset{d\in O\mathbb{Z}_{+}}{\sum }\frac{2}{d}%
x^{d}$ as follows

\begin{align*}
& \underset{d\in O\mathbb{Z}_{+}}{\sum }\frac{2}{d}x^{d} \\
& =2\left( \overset{\infty }{\underset{d=1}{\sum }}\frac{1}{d}x^{d}-\underset%
{d\in 2\mathbb{Z}_{+}}{\sum }\frac{1}{d}x^{d}\right) \\
& =2\left( \overset{\infty }{\underset{d=1}{\sum }}\frac{1}{d}x^{d}-\overset{%
\infty }{\underset{d=1\text{ }}{\sum }}\frac{1}{2d}x^{2d}\right) \\
& =2\left( -\log (1-x)+\frac{1}{2}\log (1-x^{2})\right) \\
& =\log \frac{1+x}{1-x},
\end{align*}%
where $2\mathbb{Z}_{+}=\{2,4,6,...\}$ denotes the set of all positive even
integers.

Thus we obtain

\begin{align*}
& F^{SO}(\mathcal{K};q,t;w)-F^{SO}(\mathcal{K};q,-t;w) \\
& =\underset{\mu \neq 0}{\sum }\overset{\infty }{\sum_{g=0}}\underset{\beta
\in \mathbb{Z}}{\sum }\underset{m=1}{\overset{\infty }{\sum }}\underset{k=0}{%
\overset{g}{\sum }}\frac{mn_{\mu ,g,\beta }}{\mathfrak{z}_{\mu }}\underset{%
d\in O\mathbb{Z}_{+}}{\sum }\frac{2}{d}q^{(g-2k+2m)d}t^{d\beta }pb_{d\mu }(z)
\\
& =\underset{\mu \neq 0}{\sum }\overset{\infty }{\sum_{g=0}}\underset{\beta
\in \mathbb{Z}}{\sum }\underset{m=1}{\overset{\infty }{\sum }}\underset{k=0}{%
\overset{g}{\sum }}\frac{mn_{\mu ,g,\beta }}{\mathfrak{z}_{\mu }}\underset{%
i_{1},...,i_{l(\mu )}}{\sum }\underset{d\in O\mathbb{Z}_{+}}{\sum }\frac{2}{d%
}(q^{g-2k+2m}t^{\beta }z_{i_{1}}^{\mu _{1}}z_{i_{2}}^{\mu _{2}}\cdot \cdot
\cdot z_{i_{\ell (\mu )}}^{\mu _{\ell (\mu )}})^{d} \\
& =\underset{\mu \neq 0}{\sum }\overset{\infty }{\sum_{g=0}}\underset{\beta
\in \mathbb{Z}}{\sum }\underset{m=1}{\overset{\infty }{\sum }}\underset{k=0}{%
\overset{g}{\sum }}\frac{mn_{\mu ,g,\beta }}{\mathfrak{z}_{\mu }}\underset{%
i_{1},...,i_{l(\mu )}}{\sum }\log \frac{1+q^{g-2k+2m}t^{\beta
}z_{i_{1}}^{\mu _{1}}z_{i_{2}}^{\mu _{2}}\cdot \cdot \cdot z_{i_{\ell (\mu
)}}^{\mu _{\ell (\mu )}}}{1-q^{g-2k+2m}t^{\beta }z_{i_{1}}^{\mu
_{1}}z_{i_{2}}^{\mu _{2}}\cdot \cdot \cdot z_{i_{\ell (\mu )}}^{\mu _{\ell
(\mu )}}} \\
& =\underset{\mu \neq 0}{\sum }\overset{\infty }{\sum_{g=0}}\underset{\beta
\in \mathbb{Z}}{\sum }\underset{m=1}{\overset{\infty }{\sum }}\underset{k=0}{%
\overset{g}{\sum }}\frac{mn_{\mu ,g,\beta }}{\mathfrak{z}_{\mu }}\log
\underset{i_{1},...,i_{l(\mu )}}{\prod }\frac{1+q^{g-2k+2m}t^{\beta
}z_{i_{1}}^{\mu _{1}}z_{i_{2}}^{\mu _{2}}\cdot \cdot \cdot z_{i_{\ell (\mu
)}}^{\mu _{\ell (\mu )}}}{1-q^{g-2k+2m}t^{\beta }z_{i_{1}}^{\mu
_{1}}z_{i_{2}}^{\mu _{2}}\cdot \cdot \cdot z_{i_{\ell (\mu )}}^{\mu _{\ell
(\mu )}}}
\end{align*}

Define the symmetric product as shown in the following formula:
\begin{equation*}
\big\langle1\pm \omega z^{\mu }\big\rangle=\prod_{i_{1},\ldots ,i_{\ell (\mu
)}}\Big(1\pm \omega z_{i_{1}}^{\mu _{1}}\cdots z_{i_{\ell (\mu )}}^{\mu
_{\ell (\mu )}}\Big).
\end{equation*}

Thus we have

\begin{align*}
& F^{SO}(\mathcal{K};q,t;w)-F^{SO}(\mathcal{K};q,-t;w) \\
& =\underset{\mu \neq 0}{\sum }\overset{\infty }{\sum_{g=0}}\underset{\beta
\in \mathbb{Z}}{\sum }\underset{m=1}{\overset{\infty }{\sum }}\underset{k=0}{%
\overset{g}{\sum }}\frac{mn_{\mu ,g,\beta }}{\mathfrak{z}_{\mu }}\log \frac{%
\big\langle1+q^{g-2k+2m}t^{\beta }z^{\mu }\big\rangle}{\big\langle%
1-q^{g-2k+2m}t^{\beta }z^{\mu }\big\rangle} \\
& =\log \underset{\mu \neq 0}{\prod }\overset{\infty }{\underset{g=0}{\prod }%
}\underset{\beta \in \mathbb{Z}}{\prod }\underset{m=1}{\overset{\infty }{%
\prod }}\underset{k=0}{\overset{g}{\prod }}\left( \frac{\big\langle%
1+q^{g-2k+2m}t^{\beta }z^{\mu }\big\rangle}{\big\langle1-q^{g-2k+2m}t^{\beta
}z^{\mu }\big\rangle}\right) ^{\frac{mn_{\mu ,g,\beta }}{\mathfrak{z}_{\mu }}%
}
\end{align*}

Now we obtain the infinite product formula for the orthogonal Chern-Simons
partition function of knots.

\begin{theorem}[Orthogonal Infinite Product Formula for Knots]
The Chern-Simons partition function for orthogonal quantum group invariants
can be expressed as the following infinite product formula
\begin{equation}
\frac{Z_{CS}^{SO}(\mathcal{K};q,t;w)}{Z_{CS}^{SO}(\mathcal{K};q,-t;w)}=%
\underset{\mu \neq 0}{\prod }\overset{\infty }{\underset{g=0}{\prod }}%
\underset{\beta \in \mathbb{Z}}{\prod }\underset{m=1}{\overset{\infty }{%
\prod }}\underset{k=0}{\overset{g}{\prod }}\left( \frac{\big\langle%
1+q^{g-2k+2m}t^{\beta }z^{\mu }\big\rangle}{\big\langle1-q^{g-2k+2m}t^{\beta
}z^{\mu }\big\rangle}\right) ^{\frac{mn_{\mu ,g,\beta }}{\mathfrak{z}_{\mu }}%
}
\end{equation}
\end{theorem}

\subsection{The case of a link}

Now we consider the case of link.

Given a link $\mathcal{L}$ of $L$ components, let $\overrightarrow{w}%
=(w^{1},...w^{L})$ and $\overrightarrow{z}=(z^{1},...z^{L})$ satisfying $%
w^{i}=z^{i}\ast q^{\rho }$, for $i=1,...,L$.

We generalize the symmetric product to the case of link as follows:
\begin{equation*}
\big\langle1\pm \omega (z^{1})^{\mu ^{1}}\cdot \cdot \cdot (z^{L})^{\mu ^{L}}%
\big\rangle=\prod_{i_{1,1},\ldots ,i_{1,\ell (\mu ^{1})},...,i_{L,1},\ldots
,i_{L,\ell (\mu ^{L})}}\Big(1\pm \omega \overset{L}{\prod_{\alpha =1}}\left(
(z_{i_{\alpha ,1}}^{\alpha })^{\mu _{1}^{\alpha }}\cdots (z_{i_{\alpha ,\ell
(\mu ^{\alpha })}}^{\alpha })^{\mu _{\ell (\mu ^{\alpha })}^{\alpha
}}\right) \Big),
\end{equation*}

There exist integers $n_{\overrightarrow{B};\,g,\beta }$ such that

\begin{equation*}
\overset{\infty }{\sum_{g=0}}N_{\overrightarrow{B},g,\beta }(q-q^{-1})^{g}=%
\overset{\infty }{\sum_{g=0}}n_{\overrightarrow{B},g,\beta }\underset{k=0}{%
\overset{g}{\sum }}q^{g-2k}\text{.}
\end{equation*}

In a similar way, the infinite product formula for orthogonal Chern-Simons
partition function of link can be obtained as follows

\begin{theorem}[Orthogonal Infinite Product Formula for Links]
The Chern-Simons partition function for orthogonal quantum group invariants
can be expressed as the following infinite product formula
\begin{equation}
\frac{Z_{CS}^{SO}(\mathcal{L};q,t;\overrightarrow{w})}{Z_{CS}^{SO}(\mathcal{L%
};q,-t;\overrightarrow{w})}=\underset{\overrightarrow{\mu }\neq
\overrightarrow{0}}{\prod }\overset{\infty }{\underset{g=0}{\prod }}\underset%
{\beta \in \mathbb{Z}}{\prod }\underset{m=1}{\overset{\infty }{\prod }}%
\underset{k=0}{\overset{g}{\prod }}\left( \frac{\big\langle%
1+q^{g-2k+2m}t^{\beta }(z^{1})^{\mu ^{1}}\cdot \cdot \cdot (z^{L})^{\mu ^{L}}%
\big\rangle}{\big\langle1-q^{g-2k+2m}t^{\beta }(z^{1})^{\mu ^{1}}\cdot \cdot
\cdot (z^{L})^{\mu ^{L}}\big\rangle}\right) ^{\frac{mn_{\overrightarrow{\mu }%
,g,\beta }}{\mathfrak{z}_{\overrightarrow{\mu }}}}
\end{equation}
\end{theorem}

\subsection{The case of the unknot}

In Proposition 10.2 of \cite{CC}, we have computed the free energy
associated to orthogonal Chern-Simons partition function of the knot.

\begin{equation*}
F^{SO}(\bigcirc ;q,t;z)=\underset{k=1}{\overset{\infty }{\sum }}\frac{1}{k}%
\left( 1+\frac{t^{k}-t^{-k}}{q^{k}-q^{-k}}\right) pb_{k}(z)\text{.}
\end{equation*}

Thus we have

\begin{align*}
& F^{SO}(\bigcirc ;q,t;w)-F^{SO}(\bigcirc ;q,-t;w) \\
& =\underset{k\in O\mathbb{Z}_{+}}{\sum }\frac{2}{k}\frac{t^{k}-t^{-k}}{%
q^{k}-q^{-k}}pb_{k}(w) \\
& =\underset{k\in O\mathbb{Z}_{+}}{\sum }\frac{2}{k}\frac{t^{k}-t^{-k}}{%
[k]^{2}}pb_{k}(z)
\end{align*}

Compared with (\ref{integer N}), we obtain

\begin{equation}
N_{(1),0,1}=-N_{(1),0,-1}=1.
\end{equation}
All other coefficients $N_{B;g,Q}$ are zero.

Thus we have

\begin{equation}
n_{\mu ,g,\beta }=\delta _{\mu ,(1)}\delta _{g,0}\mbox{sign}(\beta )
\end{equation}

and

\begin{equation}
\frac{Z_{CS}^{SO}(\bigcirc ;q,t;w)}{Z_{CS}^{SO}(\bigcirc ;q,-t;w)}=\underset{%
m=1}{\overset{\infty }{\prod }}\underset{i=-\infty }{\overset{\infty }{\prod
}}\left( \frac{(1+q^{2m}tz_{i})(1-q^{2m}t^{-1}z_{i})}{%
(1-q^{2m}tz_{i})(1+q^{2m}t^{-1}z_{i})}\right) ^{m}
\end{equation}

\subsection{Symmetry Property of $q\rightarrow q^{-1}$ in Infinite Product
Structure}

In this subsection, we discuss a basic symmetric property of this infinite
of product structure obtained from the orthogonal LMOV partition function.
Here we focus on the knot case only, while the case of links exactly follows
from the same analysis. In the derivation of the infinite product formula,
we assume $|q|<1 $ for the Taylor expansion of $\frac{1}{[d]^{2}}$. In the
case of $|q|>1$, the Taylor expansion is given by

\begin{equation}
\frac{1}{[d]^{2}}=\underset{m=1}{\overset{\infty }{\sum }}mq^{-2md}
\end{equation}

Therefore, the infinite product formula will be read as

\begin{eqnarray}
&&\frac{Z_{CS}^{SO}(\mathcal{K};q,t;w)}{Z_{CS}^{SO}(\mathcal{K};q,-t;w)}
\notag \\
&=&\underset{\mu \neq 0}{\prod }\overset{\infty }{\underset{g=0}{\prod }}%
\underset{\beta \in \mathbb{Z}}{\prod }\underset{m=1}{\overset{\infty }{%
\prod }}\underset{k=0}{\overset{g}{\prod }}\left( \frac{\big\langle%
1+q^{g-2k-2m}t^{\beta }z^{\mu }\big\rangle}{\big\langle1-q^{g-2k-2m}t^{\beta
}z^{\mu }\big\rangle}\right) ^{\frac{mn_{\mu ,g,\beta }}{\mathfrak{z}_{\mu }}%
}  \notag \\
&=&\underset{\mu \neq 0}{\prod }\overset{\infty }{\underset{g=0}{\prod }}%
\underset{\beta \in \mathbb{Z}}{\prod }\underset{m=1}{\overset{\infty }{%
\prod }}\underset{k=0}{\overset{g}{\prod }}\left( \frac{\big\langle%
1+q^{-g+2(g-k)-2m}t^{\beta }z^{\mu }\big\rangle}{\big\langle%
1-q^{-g+2(g-k)-2m}t^{\beta }z^{\mu }\big\rangle}\right) ^{\frac{mn_{\mu
,g,\beta }}{\mathfrak{z}_{\mu }}}  \notag \\
&=&\underset{\mu \neq 0}{\prod }\overset{\infty }{\underset{g=0}{\prod }}%
\underset{\beta \in \mathbb{Z}}{\prod }\underset{m=1}{\overset{\infty }{%
\prod }}\underset{k=0}{\overset{g}{\prod }}\left( \frac{\big\langle%
1+q^{-g+2k-2m}t^{\beta }z^{\mu }\big\rangle}{\big\langle1-q^{-g+2k-2m}t^{%
\beta }z^{\mu }\big\rangle}\right) ^{\frac{mn_{\mu ,g,\beta }}{\mathfrak{z}%
_{\mu }}}  \notag \\
&=&\underset{\mu \neq 0}{\prod }\overset{\infty }{\underset{g=0}{\prod }}%
\underset{\beta \in \mathbb{Z}}{\prod }\underset{m=1}{\overset{\infty }{%
\prod }}\underset{k=0}{\overset{g}{\prod }}\left( \frac{\big\langle%
1+(q^{-1})^{g-2k+2m}t^{\beta }z^{\mu }\big\rangle}{\big\langle%
1-(q^{-1})^{g-2k+2m}t^{\beta }z^{\mu }\big\rangle}\right) ^{\frac{mn_{\mu
,g,\beta }}{\mathfrak{z}_{\mu }}}
\end{eqnarray}

This is the symmetry of $q\rightarrow q^{-1}$ for the infinite product
formula.

\section{Acknowledgments}

We thank Pan Peng for many valuable discussions and Shengmao Zhu for giving
many helpful suggestions and for proof reading the paper.


\begin{thebibliography}{99}
\bibitem{BB} A. Beliakova, C.~Blanchet, Skein construction of idempotents in
Birman-Murakami-Wenzl algebras, \textit{Math. Ann.}, \textbf{321} (2001)
347-373.

\bibitem{BFM} V.~Bouchard, B.~Florea, M.~Marino. Topological open string
amplitudes on orientifolds, \textit{J. High Energy Phys.} \textbf{02} (2005)
002, hep-th/0411227.

\bibitem{B} R.~Brauer. On algebras which are connected with the semisimple
continuous groups, \textit{Ann. Math.} \textbf{63} (1937) 854-872.

\bibitem{Che1} L.~Chen, Chern-Simons theory of knot invariants, Ph.D thesis,
UCLA, 2009.

\bibitem{CC} L. Chen, Q. Chen, Orthogonal Quantum Group Invariants of Links,
\textit{Pacific J. of Math.}, 257 (2012), no.2, 267-318, math.QA/1007.1656

\bibitem{Che2} Q. Chen, Some mathematical aspects of quantum field theory,
Ph.D thesis, UC-Berkeley, 2009.

\bibitem{HOMFLY} P.~Freyd, D.~Yetter, J.~Hoste, W.B.R.~Lichorish, K.~Millet,
A.~Ocneanu, A new polynomial invariant of knots and links, \textit{Bull.
Amer. Math. Soc.},\textit{\ }\textbf{12} (1985) 239-246.

\bibitem{GV} R. Gopakumar and C. Vafa, On the gauge theory/geometry
correspondence, \textit{Adv. Theor. Math. Phys.}, \textbf{3} (1999)
1415-1443, hep-th/9811131.

\bibitem{J1} V. F. R.~Jones, A new knot polynomial and von Neumann algebras,
\textit{Notice Amer. Math. Soc.} \textbf{33} (1986), no. 2, 219-225.

\bibitem{J2} V. F. R. Jones, Hecke algebras representations of braid groups
and link polynomials, \textit{Ann. of Math.} \textbf{126} (1987) 335-388.

\bibitem{Kau1} L. H. Kauffman, State models and the Jones polynomial,
\textit{Topology}, \textbf{26} (1987) no. 3, 395--407.

\bibitem{Kau2} L. H.~Kauffman. An invariant of regular isotopy, \textit{%
Trans. Amer. Math. Soc.} \textbf{318} (1990) 417-471.

\bibitem{LaM} J. M. F.~Labastida, M.~Mari\~{n}o, A new point of view in the
theory of knot and link invariants, \textit{J. Knot Theory Ramif., }\textbf{%
11} (2002) 173-197.

\bibitem{LMV} J. M. F.~Labastida, M.~Mari\~{n}o, C.~Vafa. Knots, links and
branes at large N, \textit{J. High Energy Phys.}, \textbf{11} (2000) 007.

\bibitem{LiM} W.B.R.~Lichorish and K.C.~Millett, A polynomial invariant of
oriented links, \textit{Topology} \textbf{26} (1987) 107-141.

\bibitem{LLZ} C.-C. Liu, K. Liu, J. Zhou, A proof of a conjecture of Mari%
\~{n}o-Vafa on Hodge integrals, \textit{J. Differential Geom.} \textbf{65}
(2003), no. 2, 289-340.

\bibitem{LP1} K. Liu, P. Peng, On a proof of the Labastida-Mari\~{n}%
o-Ooguri-Vafa conjecture, \textit{Math. Res. Lett.}, \textbf{17} (2010), no.
3, 493-506.

\bibitem{LP2} K.~Liu, P. Peng, Proof of the Labastida-Mari\~{n}o-Ooguri-Vafa
conjecture, math-ph/0704.1526, \textit{J. Differential Geom.,} \textbf{85}
(2010) 479-525.

\bibitem{LP3} K. Liu, P. Peng, New Structure of knot invariants,
GT/1012.2636.

\bibitem{M} M. Mari\~{n}o, String Theory and the Kauffman polynomial,
\textit{Comm. Math. Phys.}, \textbf{298} (2010), no. 3, 613-643,
hep-th/0904.1088.

\bibitem{OV} H.~Ooguri, C.~Vafa, Knot invariants and topological strings,
\textit{Nucl. Phys. B.} 577 (2000) 419-438.

\bibitem{PT} J. H. Przytycki, P. Traczyk, Invariants of links of Conway
type, \textit{Kobe J. Math.}, \textbf{4} (1987) 115-139.

\bibitem{Ram1} A.~Ram, A Frobenius formula for the character of the Hecke
algebras, \textit{Invent. Math.}, \textbf{106} (1991) 461-488.

\bibitem{Ram2} A.~Ram, Characters of Brauer's Centralizer Algebras, \textit{%
Pacific J. Math.}, \textbf{169} (1995), no. 1, 173-200.

\bibitem{Ram3} A.~Ram, A "Second orthogonality relation" for characters of
Brauer algebras, \textit{European J. Combin.}, \textbf{18} (1997) 685-706.

\bibitem{RT1} N. Yu. Reshetikhin, V. G. Turaev, Invariants of 3-manifolds
via link polynomials and quantum groups. \textit{Invent. Math.} \textbf{103}
(1991), no. 3, 547--597.

\bibitem{RT2} N. Yu.~Reshetikhin, V. G.~Turaev, Ribbon graphs and their
invariants derived from quantum groups, \textit{Comm. Math. Phys}, \textbf{%
127} (1990) 1-26.

\bibitem{Wen} H.~Wenzl, On the structure of Brauer's centralizer algebra,
\textit{Ann. Math.} \textbf{128} (1988) 173-193.

\bibitem{W} E. Witten, Quantum field theory and the Jones polynomials,
\textit{Comm. Math. and Phys.}, \textbf{121} (1989) 360-379.
\end{thebibliography}
\end{document}